\documentclass[conference, 10pt]{IEEEtran}
\usepackage{HensLaTexPack2_1}
\usepackage[all]{xy}
\usepackage{bbold}
\ifCLASSINFOpdf
\else
\fi
\usepackage[cmex10]{amsmath}
\usepackage{graphicx}
\usepackage{color}
\usepackage{placeins}
\usepackage{float}
\usepackage{tabularx,colortbl}

\begin{document}
%
\title{Combinatorial Invariants of Multidimensional Topological Network Data}

\author{\IEEEauthorblockN{Gregory Henselman, Pawe\l  \, D{\l}otko}
\IEEEauthorblockA{Electrical and Systems Engineering Department\\
University of Pennsylvania\\
Philadelphia, Pennsylvania, USA\\
grh@seas.upenn.edu, dlotko@sas.upenn.edu}}


%

\newcommand{\mM}{{\mathcal M}}
\newcommand{\mX}{{\mathcal X}}
\newcommand{\mB}{{\mathcal B}}
\newcommand{\mR}{{\mathcal R}}
\newcommand{\mF}{{\mathcal F}}
\newcommand{\mD}{{\mathcal D}}
\newcommand{\mP}{{\mathcal P}}
\newcommand{\mQ}{{\mathcal Q}}
\newcommand{\mC}{{\mathcal C}}
\newcommand{\mI}{{\mathcal I}}
\newcommand{\mS}{{\mathcal S}}
\newcommand{\field}{{\mathbb F}}
\renewcommand{\d}{\partial}
\newcommand{\utf}{U^2_4}

\maketitle

\begin{abstract}
Modern applications of algebraic topology to point-cloud data analysis have motivated active investigation of combinatorial clique complexes -- multidimensional extensions of combinatorial graphs.  We show that meaningful invariants of such spaces are reflected in the combinatorial properties of an associated family of linear matroids and discuss matroid-theoretic approaches to several problems in computational topology.  Our results allow us to derive estimates of the summary statistics of  related constructs for random point cloud data, which we discuss for several sampling distributions in $\R^2$ and $\R^3$.

\end{abstract}


%
\IEEEpeerreviewmaketitle

\section{Introduction}  \label{intro}

Modern understanding of network science is broadly informed by spectral graph theory.  As early as 1952, Tutte \cite{tutteClassOfAb} demonstrated that the spectral data encoded in the node-incidence matrix of a graph can be refined through  the study of an associated chain lattice and corresponding graphic matroid.

Over the course of the past two decades, advances in applied algebraic topology have demonstrated the utility of spectral analysis of the {$\N$-graded family} of boundary operators  associated with the \emph{clique complex} of a simple graph (see \cite{EAT} and references).  Much work in this area has  focused on the descriptive statistics of  $k$th \emph{Betti numbers}, nullities of the degree-$k$ boundary operator, but the  majority of  topological data encoded in the sequence $\partial_0, \partial_1, \ld$ is lost in such abstractions.

The present paper applies matroid-theoretic methods of analysis to the study of $n$-dimensional cell complexes, and to network clique complexes in particular.  In \S\ref{theory} we discuss the topological realization of five matroid minors ($F_7, F_7^*,N_8, N_9$, and $\utf$) and four combinatorial  invariants, ($k$-connectivity, graphicness, regularity, and the max-flow min-cut property) engendered  by cellular boundary operators.  In \S\ref{numerical} we apply the results of \S\ref{theory} to obtain estimates of the statistical distribution of these minors and invariants in spaces generated by random point cloud data in $\R^2$ and $\R^3$.

A recurring theme is the observation that  topological cell complexes  efficiently encode salient features of  abstract spaces in terms of linear dependence relations, while linear matroids provide a highly developed theory for computing with such relations.  Together, these  objects offer a rich theory for modeling.

At the outset the reader should be careful to distinguish between the cellular matroid of a CW complex (definition below) and the abstract simplicial complex represented by the independence structure of a matroid.  While the latter has been  studied extensively,  we are aware of little work outside \cite{duval2012cuts} regarding the former.  Further documentation, including detailed proofs, summary statistics, and visuals may be found in the supplemental appendix \cite{multiDimNet}.

\section{Preliminaries} 

We assume familiarity with the basic definitions of cellular homology and matroid theory.  Introductory texts in both subjects abound, and  for the current discussion \cite{EAT} and \cite{MR1170126} are particularly relevant.   Throughout $M = (E, \mI)$  denotes a matroid with ground set $E$ and independence system $\mI$.  The associated rank function and circuit set are  $r$ and  $\mC$, respectively.  The dual matroid is written $M^*$, and all dual constructions are demarcated similarly.  
Linear representations are assumed to have form $\left[\; I \; | \; A \; \right]$, where $I$ is an $r(M)$-identity matrix.  It is advantageous to suppress the identity block,  placing the label of column $i = 1, \ld, r(M)$ on row $i$ of $A$.  Under this notation we say $A$ \emph{represents $M$ in standard form}.

The abstract simplex on $k+1$ vertices is denoted $\sk_k$.  The set of $k$-cells of a CW complex $X$ is denoted $X_k$.  We say $X$ is \emph{regular}\footnote{\noindent  Regularity of a CW complex should be  distinguished from  regularity of a linear matroid, an unfortunate but  unavoidable conflict of terminology.  For the purpose of this discussion we  assume  that any pair of $k$-cells in a regular CW complex share at most one $(k-1)$-dimensional face.}
 if it is compact and all attaching maps are embeddings;  it is \emph{$n$-dimensional} if the union of $n$-cells of $X$ forms an open dense subset of the ambient space.  Given $E\su X_k$, we denote the point-set topological closure in $X$ of the union of cells in $E$ by $\ov E$.  The degree-$k$  boundary operator and degree-$k$ homology group of $X$ with coefficients in $\field$ are denoted $\partial_k^\field$ and $H_k(X ; \field)$, respectively.  
 The \emph{degree-$k$ cellular matroid of $X$ with coefficients in $\field$}, denoted  $M_k ^{\field}(X)$, is the linear matroid on $\partial_k^{\field}$.  Where context leaves no room for confusion, we will drop the arguments $k$ and $\field$, and will not distinguish between $X_k$ and $M_k^{\field}(X)$.  

   Given any finite metric space $\mP \e \{p_0, \ld, p_m\}$ with metric $d$, the \emph{Vietoris-Rips complex on $\mP$ with parameter $\ek$}, denoted $\mR(\mP, \ek)$, is the clique complex on the  combinatorial graph $G = (\mP, E)$ with $E = \{ \{p_i, p_j\}: 0<d(p_i, p_j) < \ek \}$.  Where context leaves no room for confusion we will suppress the arguments $\mP$ and $\ek$.

\section{Cellular Matroids} \label{theory}

\subsection{Minors} \label{minors}

The notion of forbidden minors appears in a number of powerful classification theorems characteristic to the theory of matroids.   Because a number of computationally attractive matroid properties are invariant under reduction,\footnote{That is, deletion and contraction.} explicit construction of one or more  minors lacking a desired property is often sufficient to  obtain a negative certificate.  In special cases, a certain converse holds: a matroid is binary if and only if it has no minor isomorphic to $\utf$ \cite{tutte1958homotopy}, regular if and only if it is binary and has no minor isomorphic to $F_7$ or $F_7^*$ \cite{tutte1958homotopy, tutte1965lectures, tutte1971introduction}, and has the max-flow min-cut property with respect to a given element $l \in E$ if and only if $l$ belongs to no minor isomorphic to $F_7^*$ \cite{seymour1977matroids}.  The notation for and definition of $\utf, F_7,$ and $F_7^*$ are standard, and may be found in the appendix with those of minors $N_8$ and $N_9$, discussed below. 

Beyond classification, minor-based decomposition schemes drive a variety of test and counting algorithms with applications in combinatorics and optimization.   Each iteration of the decomposition subroutine for the min-cut algorithm of \cite{MR1170126}, for instance, is passed a minor isomorphic to $N_8$ or $N_9$, together with an explicit $k$-separation, as input.  

By a standard construction any integer matrix may be obtained as the boundary operator of a CW complex of dimension $2$ or above.  Via subdivision such constructions can be realized as simplicial complexes, so one could argue that the class of simplicial matroids is in fact almost completely general.  However, this would be to disregard the data of the underlying cellular structure.  Therefore let us restrict our attention for a moment to simplicial matroids and their minors.

Given the strong algebraic conditions imposed by the combinatorial structure of an abstract simplicial complex,\footnote{For instance, the supports of any two columns of a simplicial boundary matrix may intersect in at most one row.} one might reasonably hope that some or all of the computationally problematic minors listed above are either absent or in some sense exceptional among simplicial matroids.  However, the numerical results of \S\ref{numerical} indicate that such minors appear in abundance in random spaces, and direct computation shows that $F_7, F_7^*, N_8, N_9$ are present in $M_2^{\field_2}(\sk_5)$, and that  $\utf$ is present in $M_2^{\R}(\sk_5)$, see  appendix.  The following lemma therefore implies each minor is present in $M_l(\sk_m)$  for $l, m$ satisfying  $2 \le l \le  m - 3$.
\begin{lemma}
Any matroid realized as a minor of $M_k(\sk_n)$ can  be realized as a minor of $M_l(\sk_m)$ for $l, m$ satisfying $k \le l \le m - n \p k$.  
\end{lemma}

  It is instructive to visualize the fundamental cycles of $F_7, F_7^*, \utf$ with respect to a fixed basis in $M_2(\sk_5)$;   supporting graphics may be found in the appendix.  Lest the reader attribute special special significance to the 1-skeleton of $\sk_5$, it is shown that $F_7$ may be realized as a minor of $M_2^{\field_2}(X)$ for  $X$ a network clique complex with cliques of size 5 and below.

The following lemma is a special case of a more general result for linear matroids.  We record it here for reference in \S\ref{numerical} and \S\ref{mfmc}.

\begin{lemma} \label{bordism}
Let $X$ be an arbitrary $n$-dimensional CW complex, $Y$ be an $n$-dimensional subcomplex of $X$ with  distinguished element $y \in Y_n$, and $C$ be a circuit of $M_n(X)$ such that $C \bigcap Y_n = \{y\}$.    To every minor $N_Y$ of $M_n(Y)$ containing $y$ and every $x \in C \backslash \{y\}$ there corresponds a minor $N_X$ of $M_n(X)$ containing $x$ and a matroid isomorphism $\psi: N_X \to N_Y$ with $\psi(x) \e y$.  
\end{lemma}

Construction of the $N_X$ of Lemma \ref{bordism} should be understood as an analog of the connect sum operation in topology, and contrasted with the $k$-sum operation of \S\ref{properties}.\footnote{Here, as indicated by Lemma \ref{bordism}, circuits play the role of homotopy identity consistent with that of spheres.}   Both constructions recall a weak form of \emph{bordism}. In the dual setting, coboundary relations perform a similar function.

Since each $y \in M_2(\sk_5)$ is contained in a minor isomorphic to $F_7$ (respectively, $ F_7^*, N_8, N_9$, and $\utf$), Lemma \ref{bordism} implies that any 2-dimensional face $x$ of a Vietoris-Rips complex $\mR$ contained in a circuit that intersects any 6-clique in exactly one 2-dimensional face is also contained in a minor isomorphic to $F_7$ (respectively, $ F_7^*, N_8, N_9$, and $\utf$).  We use this fact in \S\ref{numerical} to study statistical properties of forbidden minors in random clique complexes.

\subsection{Properties} \label{properties}

\subsubsection{Connectedness}  An illustration of the conceptual utility of cellular matroids can be found in the context of $k$-connectivity.   Topological $k$-connectivity and matroid $k$-connectivity are distinct notions: a contractible CW complex is topologically $k$-connected for all $k \ge 0$, while the associated matroid fails $k$-connectivity for $k \ge 2$;  conversely the complete graph on four vertices fails topological $k$-connectivity for $k\ge 1$ while the associated matroid is 3-connected.  As it will be directly relevant to the following discussion, we provide the  definition of matroid $k$-connectivity below.  Topological $k$-connectivity is defined in \cite{MR1373655}, where the reader may also find a useful survey of connections with  the theory of greedoids.

\begin{definition}A \emph{$k$-separation} of $M \e (E, \mI)$ is a partition of $E$ into disjoint subsets $E_1, E_2$ such that $|E_1|, |E_2| \ge k$ and $r(E_1) \p r(E_2) \le  r(M) \p k-1$.   We say $M$ is \emph{$k$-connected} if $M$ has no $l$-separation for $1 \le l < k$.  
\end{definition}

The following result establishes a simple criterion for $k$-connectivity of a regular $n$-dimensional CW complex.  It should be contrasted with the homotopy theorems of Tutte \cite{tutte1958homotopy}, Maurer \cite{maurer1973matroid}, and Bj\"orner et al. \cite{bjorner1985homotopy}.  A corresponding $k$-sum operation is described in the appendix.  To simplify notation, given subsets  $E_1, E_2 \su X_{n}$, let $\nk_{n-1}(E_1, E_2)$  denote the nullity of the  map  $\tilde H_{n-1}(\ov E_1 \bigcap \ov E_2) \to \tilde H_{n-1}( \ov E_1) \oplus \tilde H_{n-1}(\ov E_2)$ induced by inclusion on relative homology groups.  
\begin{theorem} \label{rankCount}
Let $X$ be a regular $n$-dimensional CW complex and let $E_1, E_2$ be any partition of $X_n$ into disjoint nonempty subsets.  Then $r(M_n(X)) = r(E_1) \p r(E_2) \p \nk_{n-1}(E_1, E_2)$.
\end{theorem}
\begin{corollary}\label{connectedness}  For regular, $n$-dimensional cell complex $X$, the matroid $M_n(X)$ is $k$-connected if and only if $\nk_{n-1}(E_1, E_2) \ge l$ for all $1 \le l < k$ and  all partitions $X_n = E_1 \bigcup E_2$ with $|E_1|, |E_2| \ge l$.
\end{corollary}

It should be noted that if $E_1$ and $E_2$ are connected as topological spaces, then $\tilde H_0( \ov E_1) \e \tilde H_0(\ov E_2) \e 0$.  In this case $\nk_{0}(E_1, E_2)$ is the rank of $\tilde H_{0}(\ov E_1 \bigcap \ov E_2)$.   Thus, for $n\e1$, Corollary \ref{connectedness} recovers the fact that a connected simple graph is vertex $k$-connected if and only if the corresponding matroid is $k$-connected.\footnote{In the sense of \cite{MR1170126}, for $k \ge 2$, and excepting  the wheel graph $W_3$.}

Three elements of Theorem \ref{rankCount} and Corollary \ref{connectedness} bear emphasis.   First, proper interpretation of Corollary  \ref{connectedness} underscores the power of algebraic topology to reduce technical algebraic conditions to informative geometric concepts: a matroid $M_n(X)$ is $k$-connected if, in a mathematically rigorous sense, for any partition of $X_n$ into sets $E_1, E_2$ of size at least $l$ ($1 \le l < k$),  $l$ or more nontrivial loops straddling the intersection $\ov E_1 \bigcap \ov E_2$  may be contracted to points (trivial loops) in both $\ov E_1$ and $\ov E_2$.

Second, the \emph{Mayer-Vietoris sequence} of algebraic topology lies at the heart of Theorem \ref{rankCount}.   Like many constructions in this field,  Mayer-Vietoris  integrates local data distributed across spacial partitions and  dimensions.  The development of computational tools to leverage such interconnections for strictly greater information than is available via stand-alone, level-wise analysis is among the most profound contributions of algebraic topology.

Third, functoriality, a cornerstone of topological reasoning, plays a central role.  The application of functorial tools  often appears to be ``overkill'' in low-dimensional settings, however it simplifies the statements and proofs of both results listed above.

To concretize each of the remaining  matroid definitions it will be helpful to keep the following example in mind.  
\begin{definition} A \emph{flow} on matroid $M$ is a nonnegative integer-valued function $\mC \to \N$.  Given a distinguished element $l \in E$ and \emph{capacity function} $h: E \backslash\{l\} \to \N$, we say $\fk$ is \emph{admissible} if $\fk_e \le h(e)$ for all $e \in E \backslash\{l\}$.   The \emph{flow value} of $\fk$ on $e\in E$ is $\fk_e = \sum_{C \in \mC_e} \fk(C)$, where $\mC_e \e \{C \in \mC: e \in C\}$.   A flow is \emph{maximal} if its  value is maximal among admissible flows.    An \emph{$l$-cut}  is a subset $\ck \su E \backslash \{l\}$ for which $C \cap \ck\neq \emptyset$  for all $C \in \mC_l$.   The \emph{cut value} of $\ck$ with respect to $h$ is  $h(\ck) \e \sum_{e \in \ck} h(e)$.   We say $\ck$ is \emph{minimal} if $h(\ck) \le h(\ck')$ for all $l$-cuts $\ck'$.   The pair $(M,l)$ is said to have the \emph{max-flow min-cut property} (MFMC) if and only if the maximum flow value equals the minimum cut value for every capacity function $h$.
\end{definition}

A number of important combinatorial problems can be cast in terms of flows and cuts, among them the  minimum bounding surface problem, the  circuit packing problem, and the minimum cocircuit problem.

 In the context of random spaces both flows and cuts afford natural interpretations.   Posit a compact, contractible subset $\mD \su \R^2$ with connected, nonempty interior and  piecewise linear boundary, and a finite subset $\mP \su \mD$.  Heuristically, $\mD$ is a region in which one wishes to monitor the occurrence of an event, and $\mP$ is a collection of sensors distributed throughout $\mD$.  Suppose that each sensor lacks GPS capability, but may broadcast a unique identification code detectable by exactly those sensors lying within distance $\ek>0$ (a fixed parameter) and that each  may detect events occurring within distance $\ek/\sqrt{3}$.  Suppose further that a ring of fence nodes $F \su \mP$ line the perimeter of $\mD$, and that the maximal subcomplex of $\mR(\mP, \ek)$ on vertex set $F$ forms a one-dimensional circuit $\mF$ in $M_1(\mR)$.  Finally, let $sh: \mR \to \R^2$ be the linear extension of the map $\mR_0 \to \R^2$  sending $\{p_i\}$ to $p_i$;  we will call this the \emph{shadow map}.  
 
 Provided that $\tm{im}(sh|_{\mF})$ is the point-set boundary in $\R^2$ of $\mD$, \cite{MR2308949} show that $H_2(\mR, \mF) \neq 0$ implies that the image of the shadow map contains $\mD$, hence that an event anywhere in the sensing region will be detected. 
The condition $H_2(\mR, \mF) \neq 0$  holds if and only if  $H_2(\mR') \neq 0$, where $\mR'$ is a regular cell complex derived from $\mR$ by attaching a $2$-cell $l$  to $\mF$ along an embedding.  Under this construction the max-flow on $\mR'$ with respect to $l$ is equal to the maximum number of  pairwise-disjoint subsets $C_1, C_2, \ld \su \mR$ satisfying the coverage criterion, and the min-cut is the minimum number of 2-cells whose removal from $\mR$ forces $H_2(\mR') = 0$.  Insofar as one is able to measure coverage topologically, this is the tightest  lower bound on the number of triangles where sensor failure must occur for an event to go unnoticed.

\subsubsection{Graphicness} \label{graphicness}
A matroid $M$ is \emph{graphic} (respectively, cographic) if $M$ (respectively, $M^*$) can be represented by the binary node-incidence matrix of a graph.  The class of (co)graphic matroids is among the most computationally tractable, with efficient greedy algorithms for both the max-flow and min-cut problems.

We show in the appendix that any regular CW complex $X$ that embeds as the $n$-skeleton of a cellular realization of $S^{n\p1} = \{x \in \R^{n\p2}: \nn{x}_2 = 1\}$  is cographic over $\F_2$.  Thus, while the problem of determining piecewise-linear embeddability in codimension-one is NP-hard \cite{matouvsek2009hardness}, the (co)graphicness test of  Bixby and Cunningham \cite{MR594849} can provide a negative certificate for the related problem of skeletal piecewise-linear embedding in $|E|\log(|E|)$ time.  If this algorithm detects that the underlying matroid is cographic it also explicitly constructs the dual graph, which, while insufficient to prove embeddability, nevertheless provides comprehensive access to efficient graph algorithms for $M_n(X)$.

\subsubsection{Regularity}

A matroid $M$ is \emph{regular} if it can be represented by a totally unimodular (TU) matrix $A$ over $\R$.\footnote{When such is the case $A$ represents $M$ over every coefficient field.}  It is well known that integer programs with TU constraint matrices may be solved by linear relaxation;  the max-flow and min-cut problems fall under this umbrella.

It can be shown that a matroid with binary representation matrix $B$ is regular if and only if the nonzero coefficients of $B$ may be signed such that the resulting matrix is TU over $\R$.  If such a signing exists, then up to scaling of rows and columns it is unique.  There exist  efficient algorithms both to perform and to check such a signing, e.g.,\ \cite{MR1170126}.  If  this algorithm detects an incorrect signing it returns (up to permutation and scaling) a cycle matrix with 1's on the diagonal, -1's on the subdiagonal, a positive 1 in the top right entry, and zeros elsewhere.  Such matrices are of basic interest as, for example, in the case of closed cellular manifolds they represent obstructions to orientability and present as  instances of nontrivial holonomy.\footnote{In particular, such matrices are algebraic analogs of closed paths in $X$ which do not lift to closed paths in the orientable double cover, see \cite{hatcher2002algebraic}.}  Furthermore, one may deduce from the existence of such matrices that any $X$ with regular $M_k^{\field_2}(X)$ and non-TU real boundary operator operator has a real $\utf$ minor, see appendix.

It should be noted for completeness that polynomial time algorithms for considerably more general formulations of the max-flow problem are available for regular matroids, even in the absence of a linear representation \cite{hamacher1980algebraic}. 

\subsubsection{MFMC} \label{mfmc}

  The class of MFMC matroids is among the largest  of computationally tractable  isomorphism families closed under reduction, with polynomial time algorithms for both the max-flow and min-cut problems, and  substantial connections with integer optimization, e.g.\ via signed graphs \cite{MR1170126}.  By contrast, for general matroids the minimum cut and shortest cocircuit problems are NP-hard \cite{truemper1987max}. 
  
 Recall from \S\ref{minors} that a matroid fails MFMC with respect to a distinguished element $l$ if and only if $l$ is contained in a minor isomorphic to $F_7^*$.  The observations following Lemma \ref{bordism} therefore imply that $(\mR', l)$ lacks the MFMC property if there exists a circuit in $M_2^{\field_2}(\mR')$ containing $l$ and exactly one 2-dimensional face of a 6-clique of the underlying network.  We use this fact to obtain estimates on the frequency with which $(\mR', l)$ fails MFMC in \S\ref{numerical}. 
  


\section{Numerical Results}  \label{numerical}
 
The statistical distribution of algebraic invariants of randomly generated topological spaces constitutes a subject of active research, with applications in point-cloud data analysis and  notable ties to the theory of random graphs \cite{randomSimpSurvey, bubenik2010statistical}.  The simulated data described in this section are intended to serve as reference for the formulation of informed questions regarding the statistical distribution of matroid-theoretic invariants of  random space ensembles.  Point clouds are sampled from $D_1$, the uniform distribution on the unit square in $\R^2$, $D_2$, the uniform distribution on the unit square in $\R^2$ appended with fixed fence nodes satisfying the hypotheses of $F$ in the sensor  example of \S\ref{theory},  $D_3$, the standard normal distribution on $\R^2$, and $D_4$, the standard normal distribution in  $\R^3$.  For $i = 1, \ld, 4$, $n\e1, \ld, 100$, and $\ek = 0.01, \ld, 0.35$ we repeatedly sample $n$ points from $D_i$, construct the corresponding Vietoris-Rips complex with parameter $\ek$, and record the relevant summary statistics.

A brief summary of experimental results is provided below.  Further discussion, including an analysis of matroid $2$-connected components, a detailed examination of the distribution of minors within connected components, estimates on the failure of MFMC for the pair $(\mR', l)$, and detailed summary statistics and may be found in the appendix.

Our analysis is made possible by the recently published Unimodularity Library of \cite{MR3030932}.  An adaptation of this code is used to test for total unimodularity, regularity, and graphicness, and a modification of the PHAT library \cite{phat} is used to compute field homology on our own implementation of the Vietoris-Rips construction.  Further details may be found in the appendix.

\subsubsection{Connectedness}  In the majority  of samples generated by $D_2$ in which coverage holds, $M_2^{\field_2}(\mR)$ is 2-connected.  In all covered instances the number of connected components is either one or two.
\subsubsection{(Co)graphicness}  Instances of graphic and cographic matroids are absent in all distributions.  This is unsurprising given the sensitivity of such structures to minors appearing naturally in small cliques.  In light of  \S\ref{graphicness}, this provides a strong indication that complexes embeddable in $\R^3$ are rare.
\subsubsection{Regularity}  In the majority of covered instances regularity fails.  This is unsurprising for small values of parameter $\ek$ in particular, as the condition that all points in a relatively large sensing region lay within $\ek$ of a sensor but no six sensors lie pairwise within $\ek$ is quite strong.
\subsubsection{Max-flow values} Upper bounds on the maximal flow values of $(\mR', l)$  described in \S\ref{properties}  are computed by linear  relaxation of the integer program $\max \{z_l: \partial_2^\R z \e 0\}$.  Observationally, those spaces generated by dense sampling parameters\footnote{By \emph{dense} we mean those values of parameters $n$ and $\ek$ for which all collected samples engender connected matroids.}  tend to engender maximum flows of value 1.  Distributions exhibiting tenuous coverage and connectivity engender a higher percentage of max-flow values between 2 and 5.

\subsubsection{Minor $\utf$}
In the majority of cases where regularity fails, 6-cliques are present.  In rare but noteworthy instances where total unimodularity fails but regularity does not, a $\utf$ minor is present, see appendix.
\subsubsection{Minors $ F_7, F_7^*, N_8, N_9$} In the overwhelming majority of cases where regularity fails one or more 6-cliques are present, an observation consistent across distributions $D_i$.  Since each of $F_7, F_7^*, N_8,$ and $N_9$ are absent in regular complexes, it follows that the frequency with which each occurs is well approximated by that with which regularity fails.  


\section{Conclusion} \label{conclusion}
Matroid testing, classification, and decomposition theorems provide relevant information regarding related properties in topological cell complexes.  The coincidence of topological, functorial, and matroid-theoretic structures in cellular matroids suggests new questions in modeling and analysis.  In future  work we will continue to explore the descriptive  statistics of random spaces \emph{vis a vis}  forbidden minors and refined measures of connectivity.

\bibliography{matroidBib}{}

\end{document}